\documentclass[10pt]{article}
\usepackage{amsmath,amsthm,amssymb}
\providecommand{\R}{\mathbb R}

\providecommand{\E}{\mathop{\mathsf{E}{}}\nolimits}

\providecommand{\e}{\mathrm{e}}
\providecommand{\id}{\mathrm{id}}
\providecommand{\norm}[1]{\lVert #1 \rVert}
\providecommand{\drm}{\mathrm{d}}
\providecommand{\ent}{\mathop{\mathrm{H}{}}\nolimits}
\providecommand{\info}{\mathop{\mathrm{I}{}}\nolimits}
\providecommand{\shan}{\mathop{\mathrm{S}{}}\nolimits}
\providecommand{\abs}[1]{\lvert #1\rvert}
\providecommand{\sca}[1]{\langle #1\rangle}
\newtheorem{thm}{Theorem}
\newtheorem{lem}[thm]{Lemma}
\newtheorem{prop}[thm]{Proposition}
\newtheorem{cor}[thm]{Corollary}
\theoremstyle{remark}

\newtheorem*{rem}{Remark}
\theoremstyle{definition}
\newtheorem{defi}[thm]{Definition}
\title{Representation formula for the entropy and functional inequalities}
\author{
Joseph Lehec
\footnote{CEREMADE (UMR CNRS 7534) Universit\'e Paris-Dauphine.}
}
\begin{document}
\maketitle
\begin{abstract}
We prove a stochastic formula for the Gaussian relative entropy in the spirit of Borell's formula for the Laplace transform. As an application, we give simple proofs of a number of functional inequalities.
\end{abstract}
\section{Introduction: Borell's formula}
Let $\gamma_d$ be the standard Gaussian measure on $\mathbb R^d$: 
\[
\gamma_d ( \drm x) = \frac{\e^{-\abs{x}^2 /2}}{(2\pi)^{d/2}} \ \drm x 
\]
where $\abs{x} =\sqrt{x\cdot x}$ denotes the Euclidean norm of $x$. In \cite{borell1,borell2} Borell proves the following representation formula. Given a standard $d$-dimensional Brownian motion $B$ and a bounded function $f\colon \mathbb R^d\to\mathbb R$ we have
\begin{equation}
\label{borell}
\log \Bigl( \int_{\mathbb R^d} \e^f \ \drm \gamma_d \Bigr) = \sup_u \left[
 \E  \Bigl( f \bigl( B_1 + \int_0^1 u_s \ \drm s \bigr) - 
  \frac{1}{2} \int_0^1\abs{ u_s }^2 \ \drm s \Bigr) \right] ,
\end{equation}
where the supremum is taken over all random processes $u$, say bounded and 
adapted to the Brownian filtration. Among other applications, 
he derives easily the Pr\'ekopa-Leindler inequality. 
The name \emph{Borell's formula} may be unfair to
Bou\'e and Dupuis who in an earlier paper~\cite{boue-dupuis}
obtained a stronger result, allowing the function $f$ to depend
on the whole path $(B_t)_{t\in[0,1]}$
(see Theorem~\ref{boue} below for a precise statement). 
Anyway, Borell and Bou\'e-Dupuis
agree that representation formulas such as~\eqref{borell} 
arose much earlier in optimal control theory, 
particularly in Fleming and Soner's work~\cite{fleming-soner},
and Borell should definitely be credited for bringing
these techniques in the context of
functional inequalities. 

The present article deals with relative entropy. 
Let $(\Omega,\mathcal A, m)$ be a measured space 
and $\mu$ be a probability measure. 
The relative entropy of $\mu$ is defined by
\[ 
\ent ( \mu \mid m )  =
  \int_\Omega \frac{\drm \mu}{\drm m} \log\bigl( \frac{\drm \mu}{\drm m} \bigr)  \ \drm m \quad \mathrm{if}\ \mu \ll m
  \]
and $\ent( \mu \mid m ) =  +\infty$ otherwise. It is well known
that there is a Legendre duality between relative entropy and logarithmic Laplace transform:
\begin{equation}
\label{legendre}
 \ent (\mu \mid m) 
= \sup_f \Bigl( \int f \ \drm \mu  - \log\bigl( \int_\Omega \e^f \ \drm m \bigr)    \Bigr) . 
\end{equation}
The purpose of this article is to prove 
a representation formula for the Gaussian
relative entropy, both in $\mathbb R^d$ and
in the Wiener space,
providing the entropy counterparts
of the results mentioned above.
All these formulas have a common feature: Girsanov's theorem.
However, our approach is somewhat different from 
that of Borell and Bou\'e-Dupuis: 
it draws a connection with the
work of F\"ollmer~\cite{follmer1,follmer2} 
which makes the whole argument arguably simpler. 
As an application, we give new, unified and simple proofs
of a number of Gaussian inequalities. 
\section{Representation formula for the entropy}
This section contains the main results of the article. Let us recall
a couple of classical facts about relative entropy, see for instance~\cite[section 10]{varadhan}
and the references therein. If $\mathcal A$ is the Borel $\sigma$-field of
a Polish topology on $\Omega$ then it is enough to take the supremum
over bounded and continuous function in~\eqref{legendre}.
In particular the map $\mu \mapsto \ent (\mu \mid m )$ is lower semicontinuous
with respect to the topology of weak convergence of measures.
If $T \colon (\Omega , \mathcal A)\to (\Omega',\mathcal A')$ is a measurable map then 
\begin{equation}
\label{proj}
\ent ( \mu \circ T^{-1} \mid m \circ T^{-1} ) \leq \ent ( \mu \mid m ) 
\end{equation}
and assuming that $\ent ( \mu \mid m ) <+\infty$, equality occurs if and only if the density $\drm \mu / \drm m$ is a function of  $T$. 
\\
We now describe the setting of the article.
Let $\mathbb W$ be the space of continuous paths
\[
\bigl\{  w \in \mathcal C^0 (  \mathbb R_+ , \mathbb R^d ) ,  \ w_0 = 0  \bigr\} 
\]
equipped with the topology of uniform convergence on 
compact intervals. Let $\mathcal B$ be the associated Borel $\sigma$-field
and let $\gamma$ be the Wiener measure on $(\mathbb W, \mathcal B)$.
Let $x_t \colon w \mapsto w_t$ be the coordinate 
process and $(\mathcal G_t)_{t\geq 0}$ be the natural filtration of $x$. 
It is well known that $\mathcal B$ 
coincides with the smallest $\sigma$-field containing $\cup_{t\geq 0}\ \mathcal G_t$.
Let $\mathbb H$ be the Cameron-Martin space:
a path $U$ belongs to $\mathbb H$ if there exists $u\in \mathrm L^2 \bigl( [0,+\infty) ; \mathbb R^d \bigr)$
such that 
\[
U_t = \int_0^t u_s \ \drm s, \quad  t\geq 0 .
\]
The norm of $U$ in $\mathbb H$ is then defined by
\[ 
\norm{U} = \bigl(  \int_0^{+\infty} \abs{ u_s }^2 \ \drm s \bigr)^{1/2}.
\]
The Cauchy-Schwarz inequality shows that the Hilbert space $\mathbb H$
embeds continuously in $\mathbb W$. 
Given a probability space $(\Omega , \mathcal{A} , \mathsf P)$
 equipped with a filtration $(\mathcal F_t)_{t\geq 0}$ 
we call \emph{drift} any adapted
process $U$ which belongs to $\mathbb H$ almost surely.
Lastly, our Brownian motions are always $d$-dimensional, standard
and always start from $0$.
\subsection{The upper bound}
We shall use repeatedly Girsanov's formula, see~\cite[chapter 6]{liptser-shiryayev}. 
\begin{prop}
\label{majoration-prop}
Let $B$ be a Brownian motion defined on some 
filtered probability space 
$(\Omega,\mathcal A, \mathsf P , \mathcal F)$
 and let $U$ be a drift. 
Letting $\mu$ be the law of $B+U$, we have
\begin{equation}
\label{majoration-eq}
\ent ( \mu \mid \gamma ) \leq \frac{1}{2} \E \norm{U}^2 .
\end{equation}
\end{prop}
\begin{proof}
Write $U_t = \int_0^t u_s \ \drm s$ and
assume for the moment that $\norm{U}^2 = \int_0^{\infty} \abs{u_s}^2 \drm s$
is uniformly bounded. Then by Novikov's criterion
\[  M_t =  \exp\Bigl(   - \int_0^t u_s \cdot \drm B_s   -  \frac{1}{2} \int_0^t \abs{ u_s }^2 \ \drm s  \Bigr)  , \quad t\geq 0 \]
is a uniformly integrable martingale and Girsanov's formula applies. 
Under 
\[ \drm  \mathsf  Q = M_\infty \ \drm \mathsf P \]
the process $X:= B+U$ is a Brownian motion.
Therefore $X$ has law $\mu$ and $\gamma$
under $\mathsf P$ and $\mathsf Q$, respectively. 
Then by~\eqref{proj}
\[
\ent ( \mu \mid \gamma )  \leq \ent ( \mathsf P \mid \mathsf Q  ) = -  \E \log ( M_\infty) = \frac{1}{2} \E \norm{U}^2 ,
\]
which concludes the proof when $\norm{U}$ is bounded. In the general case,
define the stopping time
\[
T_n = \inf \bigl( t\geq 0 , \ \int_0^t \abs{u_s}^2 \ \drm s \geq n \bigr)  ,
\]
let $U_n$ be the stopped process $(U_n)_t = U_{t\wedge T_n}$ and $\mu_n$
be the law of $B+U_n$. 
With probability $1$
we have $\norm{U}^2 <+\infty$, thus $T_n  \rightarrow +\infty$ and 
$ U_n \rightarrow U$ in $\mathbb H$, hence in $\mathbb W$.
Therefore $\mu_n \rightarrow \mu$ weakly.
Also $\E \norm{U_n}^2 \rightarrow \E \norm{U}^2$ by monotone convergence.
Thus, using the lower semicontinuity of the entropy (observe that $\mathbb W$ is a Polish space)
\[
\begin{split}
\ent(\mu \mid \gamma) &\leq \liminf_n \ent ( \mu_n \mid \gamma) \\
 & \leq \liminf_n \frac{1}{2} \E \norm{U_n}^2  = \frac{1}{2} \E \norm{U}^2 . \qedhere
\end{split}
\]
\end{proof}
\begin{rem}
It follows immediately that when $\E \norm{U}^2 <+\infty$, 
the law of $B+U$ is absolutely continuous with respect to the Wiener measure $\gamma$. 
Let us point out that this is actually true for all drifts $U$, even if $\E \norm{U}^2 =+\infty$, see~\cite[chapter 7]{liptser-shiryayev}.
\end{rem}
\subsection{F\"ollmer's drift}
Let us address the question whether, given a probability measure $\mu$ on $\mathbb W$,
equality can be achieved in~\eqref{majoration-eq}. 
Recall that $(x_t)_{t\geq 0}$ is the coordinate process
on Wiener space $(\mathbb W, \mathcal B, \gamma)$
and that $(\mathcal G_t)_{t\geq 0}$ is its natural filtration.
The following is due to F\"ollmer~\cite{follmer1,follmer2}.
\begin{thm}
\label{follmer}
Let $\mu$ be a measure on $(\mathbb W,\mathcal B)$
having density $F$ with respect to $\gamma$.
There exists an adapted process $u$
such that under $\mu$ the following holds.
\begin{enumerate}
\item The process $U_t = \int_0^t u_s \ \drm s$ belongs to $\mathbb H$ almost surely.
\item The process $y =  x - U$ is a Brownian motion.
\item The relative entropy of $\mu$ is
\[  \ent ( \mu \mid \gamma) = \frac{1}{2}  \E^\mu \norm{U}^2 . \]
\end{enumerate}
\end{thm}
%
%
%
We sketch the proof for completeness.
\begin{proof}
Throughout $\E^\gamma$ and $\E^\mu$ denote
expectations with respect to $\gamma$ and $\mu$ respectively.
On $\mathcal G_t$ the measure $\mu$ has density
\[ F_t : = \E^\gamma \bigl(  F \ \big\vert\  \mathcal G_t \bigr)   , \]
with respect to $\gamma$. A standard martingale argument shows that
\begin{equation}
\label{truc1}
 \mu \Bigl(   \inf_{t \geq 0 }  \  F_t > 0 \Bigr)  =\mu ( F >0 ) =  1 .
 \end{equation}
Since Brownian martingales can be represented as stochastic integrals
there exists an adapted process $v$ satisfying 
\begin{equation}
\label{truc2}
\gamma \Bigl( \int_0^{+\infty} \abs{v_s}^2 \ \drm s < +\infty \Bigr) = 1
\end{equation}
and
\[
F_t  = 1 + \int_0^t  v_s  \cdot \drm x_s , \quad t\geq 0 .
\]
Let $u$ be the process defined by
\[
u_t =  \mathbf{1}_{\{F_t>0\}} \ (F_t)^{-1} v_t .
\]
It is adapted and~\eqref{truc1} and \eqref{truc2} yield
\[
\mu \Big( \int_0^\infty \abs{u_s}^2 \ \drm s < +\infty \Bigr) = 1 ,
\]
which is the first assertion of the theorem.\\
The assertion~\textit{2} follows from Girsanov's formula, see~\cite[Theorem~6.2]{liptser-shiryayev}.
 \\
Under $\mu$, we have
\[
\begin{split}
F_t & = 1 + \int_0^t F_s u_s \cdot \drm x_s \\
   & =  1 + \int_0^t F_s u_s \cdot \drm y_s  +  \int_0^t F_s \abs{u_s}^2 \ \drm s .
\end{split}
\]
Applying It\^o's formula (recall that $F$ is positive and $y$ is a Brownian motion 
under $\mu$) we obtain
\[ 
\log (F) = \int_0^{+\infty} u_s \cdot \drm  y_s +
  \frac{1}{2} \int_0^{+\infty} \abs{u_s}^2 \ \drm s . 
 \]
If $\E^\mu \norm{U}^2 < +\infty$ the local martingale part in the equation above
is integrable and has mean $0$ so that
\[ 
\ent ( \mu \mid \gamma ) = \E^\mu \log (F) =  \frac{1}{2} \E^\mu \norm{U}^2 .
\]
Again, a localization argument shows that this equality
remains valid when $ \E^\mu \norm{U}^2 = +\infty $, see~\cite[Lemma (2.6)]{follmer1}.
\end{proof}
To finish this subsection, we give a formula
for F\"ollmer's drift when the underlying density
has a Malliavin derivative, we refer to the first chapter of~\cite{nualart}
for the (little amount of) Malliavin calculus we shall use.
For suitable $F\colon \mathbb W \to \mathbb R$ we let 
 $\mathrm D F \colon \mathbb W \to \mathbb H$
 be the Malliavin derivative of $F$.
 The domain of $\mathrm D$ in the space 
 $\mathrm L^2 (\mathbb W, \mathcal B, \gamma)$
 is denoted by $\mathbb D^{2}$. 
If $F\in \mathbb D^2$ 
 then the Clark-Ocone 
formula asserts
\[ \E^\gamma (  F \mid \mathcal G_t ) =
 1 + \int_0^t \E^\gamma  \bigl( \mathrm D_s F   \mid \mathcal G_s \bigr) \cdot \drm x_s ,
 \quad t\geq 0 . \]
 We obtain the following result.
\begin{lem}
\label{malliavin}
When $F \in \mathbb D^{2}$ the process $u_t$
given by Theorem~\ref{follmer} is
\[   
u_t   =  \frac{ \E^\gamma \bigl( \mathrm D_t F \mid \mathcal G_t \bigr) }
         { \E^\gamma (  F\mid \mathcal G_t ) }  \ \mathbf{1}_{ \{ \E^\gamma (  F\mid \mathcal G_t ) >0\} }  .
\]
This implies that $\mu$-almost surely
\[
 u_t     =  \E^\mu \Bigl(  \frac{\mathrm D_t F}{F} \ \big\vert \ \mathcal G_t \Bigr) .
\]
\end{lem}
\subsection{Optimal drift in a strong sense}
According to Theorem~\ref{follmer}, 
the filtered probability space $(\mathbb W, \mathcal B, \mu, \mathcal G)$
carries a Brownian motion $y$. 
The process $x = y + U$ has law $\mu$ and the drift $U$ 
satisfies
\[
  \ent ( \mu \mid \gamma) = \frac{1}{2}  \E^\mu \norm{U}^2 .
  \]
Still, it remains open whether
\emph{given} a probability space, 
a filtration and a Brownian motion,
there exists a drift achieving equality in~\eqref{majoration-eq}. 
\\
It this section, we show that this is indeed the case,
under some restriction on the measure $\mu$.
The approach is taken from the article~\cite{baudoin} in which 
Baudoin treats the case of Brownian bridges (see subsection~\ref{bridge-section} below). 
We refer to~\cite{rogers-williams} for the background on stochastic differential
equations.
\begin{thm}
\label{yamaba}
Let $B$ be a Brownian motion defined on some 
filtered probability space $(\Omega , \mathcal A, \mathsf P , \mathcal F)$.
Let $\mu$ be a measure on $\mathbb W$,
absolutely continuous with respect to $\gamma$
and let $u_t \colon \mathbb W \to \R^d$ be the associated 
F\"ollmer process. 
If the stochastic differential equation
\begin{equation}
\label{SDE}
X_t = B_t + \int_0^t  u_s (X)  \ \drm s , \quad t\geq 0 
\end{equation}
has the pathwise uniqueness property, 
then it has a unique strong solution.
This solution $X$ satisfies the following.
\begin{enumerate}
\item The process $U_t = \int_0^t u_s(X) \ \drm s$ belongs to $\mathbb H$ almost surely.
\item The process $X$ has law $\mu$.
\item The relative entropy of $\mu$ is given by
\[  \ent ( \mu \mid \gamma ) = \frac{1}{2}  \E \norm{U}^2 . \]
\end{enumerate}
\end{thm}
\begin{proof}
According to Theorem~\ref{follmer}, on $(\mathbb W , \mathcal B, \mu)$
the coordinate process $x$ satisfies
\[
x_t = y_t + \int_0^t u_s (x)  \ \drm s 
\]
where $y$ is a Brownian motion.
Therefore~\eqref{SDE} has a weak solution.
By Yamada and Watanabe's theorem,
if pathwise uniqueness holds then~\eqref{SDE}
has a unique strong solution.
Moreover, since pathwise uniqueness
implies uniqueness in law,
the solution $X$ has law $\mu$. 
The rest of Theorem~\ref{yamaba}
concerns the law of $X$, so it is contained in 
Theorem~\ref{follmer}.
\end{proof}
We end this section by showing that
for a reasonably large class
of measures $\mu$, the stochastic differential 
equation~\eqref{SDE} does satisfy the 
pathwise uniqueness property.
\begin{defi}
\label{defS}
Let $\mathcal S$ be the class of
probability measures on $(\mathbb W, \mathcal B, \gamma)$
having a density of the form
\begin{equation}
\label{densityS}
F ( w ) = \Phi \bigl( w_{t_1} , \dotsc , w_{t_n} \bigr)
\end{equation}
for some integer $n$, for some sample $0\leq t_1 < t_2 <\dotsb < t_n$
and for some function $\Phi \colon ( \mathbb R^d )^n \to \mathbb R$ satisfying
\begin{itemize}
\item $\Phi$ is Lipschitz.
\item $\nabla \Phi$ is Lipschitz.
\item There exists $\epsilon >0$ such that $\Phi \geq \epsilon$.
\end{itemize}
\end{defi}
\begin{lem}
\label{Spath}
If $\mu$ belongs to $\mathcal S$ 
then the equation~\eqref{SDE}
has the pathwise uniqueness 
property. 
\end{lem}
\begin{proof}
Let $\mu$ have density $F$ given by~\eqref{densityS}.
Then $F\in\mathbb D^2$ and
\[
\mathrm D F (w) =  \sum_{i=1}^n  \nabla_i \Phi ( w_{t_1} , \dotsc , w_{t_n} )   \mathbf{1}_{[0,t_i]} 
\]
where $\nabla_i \Phi$ is the gradient of $\Phi$ in the $i$-th variable. 
By Lemma~\ref{malliavin}, the process associated to $\mu$ is 
\[
\begin{split}
u_t (w) & =  \frac{ \E^\gamma \bigl( \mathrm D_t F (w) \mid \mathcal G_t \bigr) }
         { \E^\gamma (  F(w) \mid \mathcal G_t ) }  \\
        &   = \sum_{i=1}^n 
          \frac{ \E^\gamma \bigl( \nabla_i \Phi ( w_{t_1} , \dotsc , w_{t_n}) \mid \mathcal G_t \bigr) }
         { \E^\gamma \bigl(   \Phi( w_{t_1} , \dotsc , w_{t_n})  \mid \mathcal G_t \bigr) }
         \mathbf{1}_{[0,t_i]}  (t) .
\end{split}     
\]
It is enough to prove that there is a constant $C$ such that
\begin{equation}
\label{lip}
 \abs{ u_t ( w ) - u_t ( \tilde{w} ) } \leq C \sup_{0\leq s\leq t} \abs{w_s - \tilde w_s} .
\end{equation}
for all $t\geq 0$ and for all $w,\tilde w \in \mathbb W$. Fix $t\geq 0$ and assume that $t_k \leq t < t_{k+1}$ for some $k\in \{0,\dotsc,n-1\}$.
By the Markov property of the Brownian motion
\[
\E \bigl( \Phi ( w_{t_1} , \dotsc , w_{t_n} ) \mid \mathcal G_t \bigr) 
 = \Psi ( w_{t_1} , \dotsc , w_{t_k} , w_t  )
\]
where $\Psi ( x_1 , \dotsc , x_{k} , x )$ equals 
\[
\int_{\mathbb W} 
\Phi \bigl( x_1 , \dotsc, x_k , x + w_{t_{k+1} - t} , \dotsc , x + w_{t_n - t} \bigr) \ \gamma ( \drm w) .
\]
Then observe that $\norm{\Psi}_{\mathrm{lip}} \leq \norm{\Phi}_{\mathrm{lip}}$. 
We have a similar property when $0\leq t < t_1$ and when $t_n\leq t$. 
The argument applies also to $\nabla_i \Phi$. The inequality~\eqref{lip} follows easily.
\end{proof}
To sum up, we have the following representation formula.
\begin{thm}
\label{representation}
Let $(\Omega,\mathcal A , \mathsf P , \mathcal F)$ be a filtered probability
space and let $B\colon \Omega \to \mathbb W$ be a Brownian motion.
For all $\mu\in\mathcal S$ we have
\[
\ent ( \mu \mid \gamma ) = \min_U \Bigl (  \frac{1}{2}  \E \norm{U}^2 \Bigr)
\]
where the minimum is on all drifts $U$ such that $B+U$ has law $\mu$. 
\end{thm}
\subsection{The Bou\'e and Dupuis formula}
In this subsection the previous results are translated in terms
of log-Laplace using the following lemma.
\begin{lem}
\label{Slegendre}
Let $f \colon \mathbb W \to \mathbb R$ bounded from below. 
For every positive $\epsilon$ there exists $\mu \in \mathcal S$ such that
\begin{equation}
\label{Seq}
\log \Bigl( \int_{\mathbb W} \e^f \ \drm \gamma \Bigr)
\leq \int_{\mathbb W} f \ \drm \mu - \ent ( \mu \mid \gamma ) + \epsilon. 
\end{equation}
\end{lem}
\begin{proof}
By monotone convergence we can assume that 
$f$ is also bounded from above, and that $\int \e^f \ \drm\gamma =1$.
Set $F = \e^f$ and let $\mu$ be a probability measure on $\mathbb W$. 
Using $t \log(t)  \leq \abs{t-1} + \abs{t-1}^2/2$ we get
\[
\begin{split}
 \ent( \mu \mid \gamma ) - \int f \ \drm \mu 
 &  \leq \int \bigl\vert \frac{G}{F} - 1 \bigr\vert  F \ \drm \gamma 
       + \frac{1}{2}  \int \bigl\vert \frac{G}{F} - 1 \bigr\vert^2 F \ \drm \gamma \\
 & \leq  \norm{F-G}_{\mathrm L^1 (\gamma) }  + C \norm{F-G}^2_{\mathrm L^2 (\gamma) }  
\end{split}
\]
where $G$ is the density of $\mu$ and $C$ is some constant (recall that $f$ is bounded below). 
Therefore, it is enough to prove that there exists $\mu \in \mathcal S$
whose density $G$ is arbitrarily close to $F$ in $\mathrm L^2(\gamma)$. 
This is left to the reader. 
\end{proof}
Here is the Bou\'e and Dupuis formula.
\begin{thm}
\label{boue}
For every function $f\colon \mathbb W \to \mathbb R$ measurable and bounded 
from below, we have
\[
 \log \Bigl( \int_{\mathbb W} \e^f \ \drm \gamma \Bigr) = \sup_U \left[
 \E  \Bigl( f ( B + U) - 
  \frac{1}{2} \norm{U}^2  \Bigr) \right] , 
   \]
where the supremum is taken over all drifts $U$.
\end{thm}
This is actually slightly more general than the result in~\cite{boue-dupuis},
which concerns the space $\mathcal {C} ( [0,T] , \mathbb R^d )$ for 
some finite time horizon $T$.
\begin{proof}
Let $U$ be a drift and $\mu$ be the law of $B+U$.
By Proposition~\ref{majoration-prop} and the entropy/log-Lapace
duality
\[
\E  \Bigl( f ( B + U) - 
  \frac{1}{2} \norm{U}^2  \Bigr)
\leq \int f \ \drm \mu - \ent ( \mu \mid \gamma ) \leq   
\log \Bigl( \int_{\mathbb W} \e^f \ \drm \gamma \Bigr) .
\]
On the other hand, given $\epsilon >0$,
there exists a probability measure
$\mu\in \mathcal S$ satisfying~\eqref{Seq}.
Since $\mu\in\mathcal S$, Theorem~\ref{representation}
asserts that there exists a drift $U$ such that
$B+U$ has law $\mu$ and satisfying
 \[ \ent ( \mu \mid \gamma ) = \frac{1}{2} \E \norm{U}^2  . \]
Then~\eqref{Seq} becomes
\[
\log \Bigl( \int_{\mathbb W} \e^f \ \drm \gamma \Bigr) \leq
 \E  \Bigl( f ( B + U) -  \frac{1}{2} \norm{U}^2  \Bigr) + \epsilon ,
 \]
which concludes the proof.
\end{proof}
\subsection{Brownian bridges}
\label{bridge-section}
A measure $\mu$ on $\mathbb W$ satisfying
\begin{equation}
\label{bridge-def}
\mu ( \drm w ) =  \rho(w_1) \ \gamma ( \drm w)
\end{equation}
where $\rho$ is some density on $(\mathbb R^d , \gamma_d)$
is said to be a Brownian bridge. It can be seen as the law of a Brownian
motion conditioned to have law 
$\rho (x) \gamma_d (\drm x)$
at time $1$. 
\begin{lem}
\label{bridge-lemma}
Let $\nu$ have density $\rho$ with respect to $\gamma_d$, we have
\[ \ent ( \nu \mid \gamma_d ) = \inf_{\mu}  \Big( \ent ( \mu \mid \gamma ) \Bigr) \]
where the infimum is on all probability measures satisfying $\mu \circ  ( x_1 ) ^{-1} = \nu$. 
The infimum is attained when $\mu$ is the bridge~\eqref{bridge-def}.
\end{lem}
In other words, among all processes having law $\nu$
at time $1$, the bridge minimizes the relative entropy.
This is essentially a particular case of~\eqref{proj}, 
see also~\cite{baudoin} and~\cite[page 161]{follmer3}.\\
Assume that $\rho$ is differentiable and that $\nabla \rho \in \mathrm L^2 ( \gamma_d)$.
Then $F (w) = \rho (w_1)$ belongs to $\mathbb D^2$ and has Malliavin derivative
\[
\mathrm D F (w) = \nabla \rho ( w_1) \mathbf{1}_{[0,1]} .
\]
By Lemma~\ref{malliavin} the F\"ollmer process of the bridge $\mu$ is such that
\[
 u_t = \E^\mu \bigl( \nabla \log(\rho) (w_1 ) \mid \mathcal G_t \bigr) \mathbf{1}_{[0,1]} (t)  , \quad \mu-a.s. 
\]
We obtain the following result.
\begin{lem}
\label{martingale}
Under $\mu$, the process $(u_t)_{t\in [0,1]}$ is a martingale. In particular
\[
\E^\mu ( u_t ) = \E^\mu \nabla  \log(\rho)  (w_1 ) = 
\E^\gamma \nabla \rho (w_1) = \int_{\mathbb R^d} x \ \nu (\drm x). 
\]
\end{lem}
Now  assume that $\rho$ and $\nabla \rho$ are Lipschitz
and that $\rho \geq \epsilon$, so that the bridge $\mu$ belongs to $\mathcal S$.
It is easily seen that $u_t$ can also be written as
\[ 
u_t (w) = \nabla \log P_{1-t} \rho ( w_t )   \ \mathbf{1}_{[0,1]} (t) ,
\]
where $P_t$ denotes the heat semigroup
on $\mathbb R^d$:
\[  \partial_t P_t = \frac{1}{2} \Delta P_t . \]
The stochastic differential equation~\eqref{SDE} becomes
\begin{equation}
\label{SDEbis}
 X_t = B_t + \int_0^{t\wedge 1} \nabla \log ( P_{1-s}  \rho ) ( X_s) \ \drm s , \quad t\geq 0 .
\end{equation}
By Lemma~\ref{Spath}, there is a unique strong solution. 
Combining Lemma~\ref{bridge-lemma} with Theorem~\ref{yamaba}
we obtain the following dual formulation
of Borell's result~\eqref{borell}.
\begin{thm}
\label{bridge-thm}
Let $\nu$ and $\rho$ be as above. Then
\[ \ent ( \nu \mid \gamma_d ) = \inf_U \Bigl( \frac{1}{2} \E \norm{U}^2 \Bigr) \]
where the infimum is taken on all drifts $U$ satisfying $B_1 + U_1 = \nu$ in law.
The infimum is attained by the drift
\[ 
U_t = \int_0^{t\wedge 1} \nabla \log ( P_{1-s}  \rho ) ( X_s) \ \drm s , 
\]
where $X$ is the unique solution of~\eqref{SDEbis}.
\end{thm}
%
%
%
\section{Applications}
Following Borell, we now derive functional 
inequalities from the representation formula.
Let us point out that in all but one applications 
we use Proposition~\ref{majoration-prop} and Theorem~\ref{follmer}  rather
than Theorem~\ref{representation}.
\subsection{Transportation cost inequality}
Let $\mathrm{T}_2$ be the transportation cost
for the Euclidean distance squared: given
two probability measures $\mu$ and $\nu$ on 
$\mathbb R^d$
\begin{equation}
\label{transport-def}
\mathrm{T}_2 ( \mu, \nu ) = \inf \Bigl( \int_{\mathbb R^d \times \mathbb R^d}  \abs{x-y}^2 \ \drm \pi(x,y) \Bigr)^{1/2} ,
\end{equation}
where the infimum is taken over all couplings $\pi$ of $\mu$ and $\nu$,
namely all probability measures on the product space $\mathbb R^d \times \mathbb R^d$ 
having marginals $\mu$ and $\nu$. 
There is a huge literature about this optimization problem,
usually referred to as Monge-Kantorovitch problem, see Villani's book \cite{villani}. 
Talagrand's inequality asserts that 
\[ 
\mathrm{T}_2 (  \nu , \gamma_d) \leq 2 \ent( \nu \mid \gamma_d ) 
\]
for every probability measure $\nu$ on $\mathbb R^d$. 
The purpose of this subsection is to prove a Wiener space version of this inequality.

On Wiener space the natural definition of $\mathrm T_2$
involves the norm of the Cameron-Martin space $\mathbb H$:
given two probability measures $\mu,\nu$ on $(\mathbb W, \mathcal B)$
\[
\mathrm{T}_2 ( \mu, \nu ) = \inf \Bigl(
\int_{\mathbb W \times \mathbb W}  \norm{w - w'}^2 \ 
   \pi( \drm w , \drm w') \Bigr) ,
\]
where the infimum is taken over all couplings $\pi$ of 
$\mu$ and $\nu$ such that $w-w' \in \mathbb H $ for $\pi$-almost all $(w,w')$.
\begin{thm}
Let $\mu$ be a probability measure on
$(\mathbb W,\mathcal B)$. Then
\[ \mathrm{T}_2 (  \mu , \gamma ) \leq 2 \ent( \mu \mid \gamma )  . \]
\end{thm}
Here is a short proof based of Theorem~\ref{follmer}. 
Fair enough, Feyel and {\"U}st{\"u}nel~\cite{feyel-ustunel}
have a very similar argument.
\begin{proof}
Assume that $\mu$ is absolutely 
continuous with respect to $\gamma$ 
(otherwise $\ent(\mu \mid \gamma)=+\infty$).
According to Theorem~\ref{follmer}
there exists a Brownian motion 
$B$ and a drift $U$ such that $B+U$ has law $\mu$ and
\[  \ent ( \mu \mid \gamma ) =  \frac{1}{2} \E \norm{U}^2  . \]
Then $(B,B+U)$ is a coupling of $(\gamma,\mu)$
and by definition of $\mathrm{T}_2$
\[
 \mathrm{T}_2 (\mu,\gamma)^2  \leq \E \norm{ U }^2  
    = 2 \ent ( \mu \mid  \gamma) . \qedhere
\]
\end{proof}
Let us point out that
Talagrand's inequality 
can be recovered
easily from this theorem, applying it
to a Brownian bridge. Details are left to the reader.
\subsection{Logarithmic Sobolev inequality}
In this section we prove
the logarithmic Sobolev inequality 
for the Wiener measure,
which extends the classical 
log-Sobolev inequality 
for the Gaussian measure, due to Gross~\cite{gross}.
When $\mu$ is a measure on $(\mathbb W,\mathcal B,\gamma)$ with density 
$F$ such that $\mathrm D F$ is well defined,
the Fisher information of $\mu$ is
\[ \info(\mu \mid \gamma ) = \int_{\mathbb W} \frac{\norm{\mathrm D F}^2}{F} \ \drm \gamma
                      = \int_{\mathbb W}  \Bigl\Vert \frac{\mathrm D F}{F} \Bigr\Vert^2 \ \drm \mu . \]
\begin{thm}
Let $\mu$ have density $F$ with respect to $\gamma$
and assume that $F \in\mathbb D^2$. Then
\begin{equation}
\label{log-sob-wiener}
 \ent (\mu \mid \gamma )
 \leq \frac{1}{2} \info (\mu \mid \gamma) .
 \end{equation}
\end{thm}
\begin{proof}
We consider the probability space $(\mathbb W, \mathcal B,\mu)$.
Recall that $(\mathcal G_t)_{t\geq 0}$ is the filtration of the coordinate process. 
By Theorem~\ref{follmer} and Lemma~\ref{malliavin}, letting
\[ u_t =   \E^\mu \Bigl( \frac{\mathrm D_t F}{F} \mid \mathcal G_t  \Bigr) \]
we have 
\[
\ent ( \mu \mid \gamma ) = 
\frac{1}{2} \E^\mu  \int_0^\infty \abs{ u_t }^2 \ \drm t .
\]
By Jensen's inequality
\[ \E^\mu \abs{ u_t  }^2 \leq  \E^\mu  \Bigl\vert \frac{\mathrm D_t  F}{F} \Bigr\vert^2  \]
so that
\[
\ent ( \mu \mid \gamma ) \leq 
\frac{1}{2} \E^\mu \Bigl\Vert \frac{\mathrm D F}{F} \Bigr\Vert^2  
\]
which is the result.
\end{proof}
This may not be the most straightforward 
proof, see~\cite{capitaine}. 
Let us emphasize that applying~\eqref{log-sob-wiener} 
to a Brownian bridge yields the usual log-Sobolev inequality. 
More precisely, let $\nu$ be a probability measure on 
$\mathbb R^d$ having a smooth density $\rho$ with respect to $\gamma_d$ 
and let $\mu$ be the measure on $\mathbb W$ given by
\[  \mu ( \drm w )  = \rho ( w_1 ) \ \gamma( \drm w ) . \]
Then $\ent ( \nu \mid \gamma_d ) = \ent ( \mu \mid \gamma )$. On the other hand
letting $F ( w ) = \rho (w_1 )$ we have
\[ \mathrm D F (w) = \nabla \rho ( w_1 ) \ \mathbf{1}_{[0,1] } , \]
which implies easily that $\info( \nu \mid \gamma_d ) = \info (\mu \mid \gamma ) $.
Thus~\eqref{log-sob-wiener} becomes
\[
\ent ( \nu \mid \gamma_d ) \leq
\frac{1}{2} \info( \nu \mid \gamma_d ) .
\]
\subsection{Shannon's inequality}
Given a random vector $\eta$ on $\mathbb R^d$
having density $\rho$ with respect to the Lebesgue 
measure, Shannon's entropy is defined as
\[
\shan ( \eta ) = - \int_{\mathbb R^d} \rho \log( \rho ) \ \drm x .
\]
In other words $\shan ( \eta ) = - \ent ( \nu \mid \lambda_d )$ where $\nu$ is the law of $\eta$ and $\lambda_d$
is the Lebesgue measure on $\mathbb R^d$.
\begin{thm}
Let $\eta,\xi$ be independent random vectors on $\R^d$ and $\theta \in [0,\pi/2]$
\begin{equation} 
\label{shannon}
\shan ( \cos(\theta) \eta + \sin(\theta) \xi ) \geq \cos( \theta )^2 \shan ( \eta ) +\sin (\theta)^2  \shan ( \xi ).
\end{equation}
\end{thm}
This inequality plays a central role in information theory, 
see~\cite{dembo-al} for an overview on the topic. 
\begin{proof}
Let $\nu_\theta$ be the law of $\cos(\theta) \eta + \sin(\theta) \xi$.  
By Theorem~\ref{follmer}, Lemma~\ref{bridge-lemma} and Lemma~\ref{martingale} there exists
a Brownian motion $X$ and a drift $U$ such that
\begin{itemize}
\item $X_1 + U_1$ has law $\nu_0$.
\item $\ent ( \nu_0 \mid \gamma_d ) =  \E \norm{U}^2 /2$.
\item $\E (U) = \E ( \eta )  \ \mathbf{1}_{[0,1]}$.
\end{itemize}
Similarly, there exists a Brownian motion $Y$ and a drift $V$ 
satisfying the corresponding properties for $\nu_{\pi/2}$. 
Besides, we can clearly assume that $Y$ is independent of $X$.
Then $\cos(\theta) X +\sin(\theta) Y$ is a Brownian motion
and
\[
\cos(\theta) X_1  +\sin(\theta) Y_1 + \cos(\theta) U_1 +\sin(\theta) V_1
\]
has law $\nu_\theta$. By Proposition~\ref{majoration-prop} and Lemma~\ref{bridge-lemma}
\[
 \ent ( \nu_\theta \mid \gamma_d )  \leq \frac{1}{2} \E \norm{ \cos(\theta) U +\sin(\theta) V }^2  .
 \]
Denoting the inner product in $\mathbb H$ by $\sca{\cdot,\cdot}$ we have 
\[\E \sca{U,V} = \sca{\E U , \E V} = ( \E \eta ) \cdot ( \E \xi ) , \]
so that
\[
 \begin{split}
  \ent ( \nu_\theta \mid \gamma_d ) \leq  \cos( \theta )^2 \ent ( \nu_0 \mid \gamma_d )  
   & +  \sin(\theta)^2 \ent( \nu_{\pi/2} \mid \gamma_d ) \\
  &   + \cos(\theta) \sin(\theta) \ ( \E \eta ) \cdot ( \E \xi ) .
\end{split}
\]
This is easily seen to be equivalent to~\eqref{shannon}.
\end{proof}
\subsection{Brascamp-Lieb inequality}
Let us focus on a family of inequalities dating back to Brascamp and Lieb's
article~\cite{brascamp-lieb} on optimal constants in Young's inequality. 
Since then a number of nice alternate proofs have been discovered, 
see~\cite{barthe-huet,carlen-cordero} and the survey article~\cite{ball-handbook}. 
This subsection is inspired by the (unpublished) proof of Maurey relying on Borell's formula.

Let $E$ be a Euclidean space, let $E_1,\dotsc,E_m$ be subspaces and for all $i$ let $P_i$ be the orthogonal projection with range $E_i$. The crucial hypothesis is the so-called frame condition: there exist $c_1,\dots,c_m$ in $\mathbb R_+$ such that
\begin{equation}
\label{frame}
\sum_{i=1}^m c_i P_i = \id_E.
\end{equation}
Let $x\in E$, we then have $\abs{x}^2 = \bigl(\sum c_i P_i x \bigr) \cdot  x$ and since $P_i$ is an orthogonal projection
\begin{equation}
\label{frame1}
\abs{x}^2 =\sum_{i=1}^m c_i \abs{P_i x}^2 .
\end{equation}
From now on $\mathbb W$ denotes the space of continuous paths taking values in $E$ and starting from $0$ and $\gamma$ denotes the Wiener measure on $\mathbb W$. The spaces $\mathbb W_i$ and measures $\gamma_i$ are defined similarly.
\begin{thm}
\label{BLH}
Under the frame condition, for every probability measure $\mu$ on $\mathbb W$ we have
\[ \ent( \mu \mid \gamma ) \geq \sum_{i=1}^m c_i \ent ( \mu_i \mid \gamma_i )  , \]
where $\mu_i = \mu \circ P_i^{-1}$ is the push-forward of $\mu$ by the projection $P_i$. 
\end{thm}
\begin{proof}
According to Theorem~\ref{follmer} there exists
a standard Brownian motion $B$ on $E$
and a drift $U$ such that $B+U$ has law $\mu$ and
\[
\ent ( \mu \mid \gamma ) = \frac{1}{2} \E\norm{U}^2 . 
\]
Since $P_i$ is an orthogonal projection, the process $P_i B$
is a standard Brownian motion on $E_i$. Also $P_i  B  + P_i U$ has law $\mu \circ P_i^{-1} = \mu_i$. 
By Proposition~\ref{majoration-prop} 
\[ \ent ( \mu_i \mid \gamma_i ) \leq \frac{1}{2}  \E \norm{P_i U}^2 , \quad i = 1 , \dotsc, m.\]
On the other hand, the frame condition~\eqref{frame1} implies easily that
\[
\norm{U}^2 =  \sum_{i=1}^n  c_i \norm{P_i U}^2 
\]
pointwise. Taking expectation yields the result.
\end{proof}
As observed by Carlen and Cordero~\cite{carlen-cordero}, this super-additivity property of the relative entropy is equivalent to the following Brascamp-Lieb inequality.
\begin{cor}
\label{BL}
Under the frame condition, given $m$ functions $F_i\colon \mathbb W_i  \to \R_+$, we have
\[  \int_{\mathbb W} \prod_{i=1}^m (  F_i\circ P_i )^{c_i}  \ \drm \gamma  \leq
\prod_{i=1}^m \Bigl( \int_{\mathbb W_i} F_i \ \drm \gamma_i \Bigr)^{c_i}  . \]
\end{cor}
When the functions $F_i$ depend only on 
the point $w_1$ rather than on the whole path
$w$ we recover the usual Brascamp-Lieb inequality
for the Gaussian measure.
\subsection{Reversed Brascamp-Lieb inequality}
Again $E$ is a Euclidean space and $E_1,\dotsc,E_m$ are subspaces satisfying the frame condition~\eqref{frame}. Observe that if $x_1,\dotsc,x_m$ belong to $E_1 ,\dotsc,E_m$ respectively, then for any $y\in E$, the Cauchy-Schwarz inequality and \eqref{frame1} yield
\[ 
\begin{split}
 \bigl( \sum_{i=1}^m c_i x_i \bigr) \cdot  y 
 & = \sum_{i=1}^m c_i  ( x_i  \cdot  P_i y )  \\
 & \leq 
\Bigl( \sum_{i=1}^m c_i \abs{ x_i }^2 \Bigr)^{1/2} 
\Bigl( \sum_{i=1}^m c_i \abs{ P_i y }^2 \Bigr)^{1/2} \\
& = \Bigl( \sum_{i=1}^m c_i \abs{ x_i }^2 \Bigr)^{1/2} \abs{y} .
\end{split}
\]
Hence
\begin{equation}
\label{frame2}
\Bigl\vert \sum_{i=1}^m c_i x_i \Bigr\vert^2 \leq 
 \sum_{i=1}^m c_i \abs{ x_i }^2  .
\end{equation}
Let $\mathcal S_i$ be the class of probability
measures on $E_i$ which satisfy the conditions of
Definition~\ref{defS}, replacing $\mathbb R^d$ by $E_i$. 
Here is the reversed version of Theorem~\ref{BLH}.
\begin{thm}
\label{RBLH}
Given $m$ probability measures $\mu_1,\dotsc,\mu_m$ belonging to
$\mathcal{S}_1,\dotsc,\mathcal{S}_m$ respectively, there exist $m$ processes $X_1,\dotsc,X_m$
(defined on the same probability space) such that
\begin{enumerate}
\item $X_i$ has law $\mu_i$ for all $i=1 ,\dotsc, m$.
\item Letting $\mu$ be the law of $\sum c_i X_i$ we have
\[
\ent \bigl( \mu \mid \gamma ) \leq 
\sum_{i=1}^m c_i \ent ( \mu_i \mid \gamma_i ) .
\]
\end{enumerate}
\end{thm}
\begin{proof}
Again let $B$ be a standard Brownian motion on $E$. 
For $i=1,\dotsc,m$, the process $P_i B$ is a standard Brownian motion on $E_i$.
Since $\mu_i \in \mathcal S_i$ there exists a drift $U_i$ such that the process
$X_i  = P_i B +  U_i$ has law $\mu_i$ and
\[  \ent ( \mu_i \mid \gamma_i) = \frac{1}{2} \E \norm{U_i}^2 . \] 
Let $X = \sum c_i X_i$ and let $\mu$ be the law of $X$. Since $\sum c_i P_i$ is the identity of $E$
\[  X = B + \sum_{i=1}^m c_i U_i .   \]
By Proposition~\ref{majoration-prop}, we get
\[ 
\ent ( \mu \mid \gamma )  \leq \frac{1}{2} \E \bigl\lVert \sum_{i=1}^m c_i U_i \bigr\rVert^2 . 
\]
On the other hand~\eqref{frame2} easily implies that
\[  \bigl\lVert \sum_{i=1}^m c_i U_i \bigr\rVert^2 \leq  \sum_{i=1}^m c_i \norm{U_i}^2 , \]
pointwise. Taking expectation we get the result. 
\end{proof}
This sub-additivity property of the entropy is a multi-marginal
version of the \emph{displacement convexity} property put 
forward by Sturm~\cite{sturm}. 
By duality, we obtain the following reversed Brascamp-Lieb inequality.
\begin{cor}
\label{RBL}
Assuming the frame condition, given $m$ functions $F_i \colon \mathbb W_i\to \R_+$ bounded away from $0$,
and a function $G\colon \mathbb W \to \R_+$ satisfying
\begin{equation}
\label{hyprbl}
 \prod_{i=1}^m   F_i(w_i)^{c_i} \leq G \Bigl( \sum_{i=1}^m c_i w_i \Bigr) 
\end{equation}
for all $(w_1,\dotsc,w_m)\in \mathbb W_1\times\dotsb\times \mathbb W_m$, we have 
\[ \prod_{i=1}^m  \Bigl( \int_{\mathbb W_i}  F_i \ \drm \gamma_i   \Bigr)^{c_i}
   \leq \int_{\mathbb W} G \ \drm \gamma . \]
\end{cor}
\begin{proof}
By Lemma~\ref{Slegendre}, for every $i$, there exists a measure $\mu_i \in\mathcal S_i$ such that
\[ 
\log \Bigl( \int_{\mathbb W_i} F_i \ \drm \gamma_i \Bigr)
\leq \int_{\mathbb W_i} \log(F_i) \ \drm \mu_i -  \ent ( \mu_i\mid \gamma_i ) +\epsilon .
\] 
Let $X_1,\dotsc,X_m$ be the random processes given by the previous theorem, 
let  $X = \sum c_i X_i$ and let $\mu$ be the law of $X$. 
Then by duality and the hypothesis~\eqref{hyprbl} we get
\[
\begin{split}
\log\Bigl( \int_{\mathbb W} G \ \drm \gamma \Bigr) & \geq \E \log(G) (X) - \ent ( \mu \mid \gamma ) \\
 & \geq \E \bigl( \sum_{i=1}^m c_i \log(F_i) (X_i) \bigr) - \ent( \mu \mid \gamma ) .
\end{split}
\]
Since $\ent( \mu \mid \gamma ) \leq \sum c_i \ent (\mu_i \mid \gamma_i )$, this is at least
\[
\sum c_i  \Bigl( \log \bigl( \int_{\mathbb W_i} F_i \ \drm \gamma_i  \bigr) - \epsilon \Bigr) .
\]
Letting $\epsilon$ tend to $0$ yields the result. 
\end{proof}
Again when the functions depend only on the value of the path at time $1$, we
recover the reversed Brascamp-Lieb inequality for the Gaussian measure,
which is due to Barthe~\cite{barthe-invent}.
%
%
\begin{center}
\textbf {Acknowledgements}
\end{center}
The author is grateful to Patrick Cattiaux and  Massimiliano Gubinelli 
for communicating references and
to Christian L\'eonard, Bernard Maurey and Patrick Cattiaux again for valuable discussions.

\end{document}